
\input amstex

\input epsf.tex

\documentstyle{amsppt}

\overfullrule=0pt

\define\sset{\subseteq}
\define\ctln{\centerline}

\define\blkbox{{\vrule height1.5mm width1.5mm depth1mm}}

\topmatter

\title\nofrills Families of knots for which Morton's inequality is strict\endtitle
\author  Mark Brittenham$ ^1$ and Jacqueline Jensen$ ^2$ \endauthor

\leftheadtext\nofrills{Mark Brittenham and Jacqueline Jensen}
\rightheadtext\nofrills{Knots with strict Morton's inequality}

\address   Department of Mathematics, University of Nebraska, 203 Avery Hall, Lincoln, NE 68588-0130  \endaddress
\email   mbrittenham2\@math.unl.edu \endemail

\address   Department of Mathematics and Statistics, Sam Houston State University, Huntsville, Texas 77341  \endaddress
\email jensen\@shsu.edu \endemail


\thanks  \hskip1pc$ ^{1}$ Research supported in part by NSF grant \# DMS$-$0306506 \endthanks
\thanks  \hskip1pc$ ^{2}$ Research supported in part by NSF grant \# DMS$-$0354281 \endthanks


\abstract 
We describe a procedure for creating infinite families of knots,
each having the maximum degree of their HOMFLY polynomial strictly
less than twice their canonical genus.
\endabstract

\endtopmatter

\document

\heading{\S 0 \\ Introduction}\endheading

Every knot $K$ in the 3-sphere $S^3$ is the boundary of a compact
orientable surface $\Sigma\sset S^3$, known as a Seifert surface 
for the knot $K$. In fact, in the 1930's Seifert [Se] gave an 
algorithm which, given a diagram of the knot $K$, produces
such a surface. The algorithm consists of orienting the knot
diagram, breaking each crossing and reconnecting the four ends
according to the orientation, without re-introducing a crossing, 
producing disjoint ``Seifert'' circles in the projection plane, bounding (after
offsetting nested circles) disjoint ``Seifert'' disks, and then reintroducing
the crossings by stitching the disks together with half-twisted
bands.

Every knot therefore has infinitely many such
``canonical'' surfaces (of arbitrarily high genus). The minimum of 
the genera of the surfaces built by Seifert's algorithm, over 
all diagrams of the knot $K$, is known as the {\it canonical genus} 
or {\it diagrammatic genus} of $K$, denoted $g_c(K)$. The minimum
genus over all Seifert surfaces, whether built by Seifert's
algorithm or not, is known as the {\it genus} of $K$, and
denoted $g(K)$.

In 1985, shortly after the discovery of the HOMFLY polynomial [FHLMOY],
Morton [Mo] showed that the highest degree of one of the two
variables gave a lower bound on $2g_c(K)$; details are outlined
below. This was perhaps the first piece of information encoded in 
the HOMFLY polynomial to be related to topological information 
about the knot $K$. Morton's inequality, $M(K) = $maxdeg$_z P_K(v,z)\leq 2g_c(K)$ ,
has since been shown to be an equality for many classes of knots.
These include all of the knots having 12 or fewer crossings [St2],
all alternating knots [Cr],[Mu], and, more generally, all homogeneous knots [Cm],
and the Whitehead doubles of 2-bridge knots [Na],[Tr] and pretzel knots [BJ].
It wasn't until 1998 and later that knots were found for which Morton's
inequality was strict.
The first such were found by Stoimenow [St2],[St3] while analyzing 
the survey of knots through 16 crossings. 

In this paper we show how to use Stoimenow's second collection of examples, or any other
example that might be built along the same lines, to build
infinite families of knots having $M(K) = $maxdeg$_z P_K(v,z) < 2g_c(K)$.
In particular, we have the following

\proclaim{Theorem}
Suppose $K$ is a knot with $g(K)=g_c(K)$ and $M(K) < 2g_c(K)$.
Let $D$ be a $g_c$-minimizing diagram of $K$, and suppose there is a 
crossing $c$ in $D$ which bounds a half-twisted band connecting
a pair of Seifert disks, and the knot $K^\prime$ obtained by changing
the crossing $c$ has $g_c(K^\prime) < g_c(K)$. Then the knots $K_n$,
bounding the canonical Seifert surfaces $\Sigma_n$
obtained by replacing the half-twisted band at the crossing $c$ with
$2n+1$ half-twisted bands in parallel, all joining the same pair of 
Seifert disks, all satisfy $M(K_n) < 2g_c(K_n) = 2(g_c(K)+n)$.
\endproclaim

The figures in section 2 should make the construction clear, if the description
above did not. We note that the condition $g_c(K^\prime) < g_c(K)$ simply states
that the diagram $D^\prime$ obtained from changing the crossing
$c$ is not a $g_c$-minimizing diagram for $K^\prime$; this can often in 
practice be verified fairly quickly, as our examples below show.

In particular, since, as we shall see, some of Stoimenow's examples satisfy 
the hypotheses 
of the theorem, and we only need one to get the process started, we have

\proclaim{Corollary} There exist knots $K$ with arbitrarily large
canonical genus for which $M(K) < 2g_c(K)$ ;
in particular, there are infinitely many of them.
\endproclaim

\heading{\S 1 \\ Notations and preliminaries}\endheading

$K$ will always denote a knot in the 3-sphere $S^3$, 
$N(K)$ a tubular neighborhood of $K$, $E(K) = S^3 \setminus$ int$N(K)$ the
exterior of $K$, $\Sigma$ a Seifert surface for $K$ (which we treat as embedded
in $S^3$, or properly embedded in $E(K)$, as needed), and $E(\Sigma) = 
E(K)|\Sigma = E(K)\setminus$
int$N(\Sigma)$ the exterior of $\Sigma$ .
$E(\Sigma)$ can be thought of as a {\it sutured manifold}
[Ga1], that is, a compact manifold $M$ with $\partial M = 
R_+\cup R_-$, where $R_+\cap R_- = \gamma$ is a collection of simple 
closed curves, the {\it sutures}. Formally, torus boundary components 
with no sutures may also be treated as a suture.
In our case $\gamma$ = the core of the annulus
$\partial E(K)\setminus \partial\Sigma$, cutting $\partial E(\Sigma)$
into two copies of $\Sigma$, pushed off of $\Sigma$ to the $+$- and $-$-sides.
The pair $(E(K),\partial E(K))$ is also an example.

\leavevmode

\epsfxsize=4in
\centerline{{\epsfbox{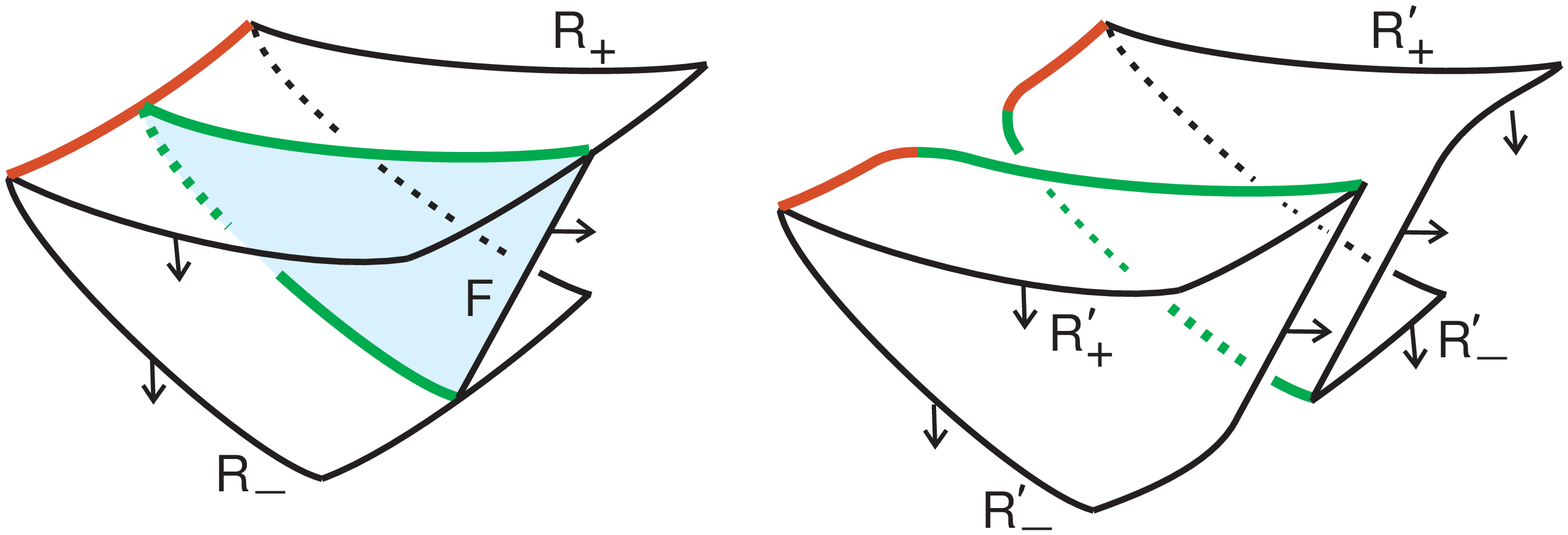}}}

\centerline{Figure 1: Sutured manifold decomposition}

\medskip

The theory of sutured manifolds will play a central role to our 
proofs below; we summarize the main points here. $R_+$ and $R_-$
are thought of as having normal orientations, pointing into $M$ along
$R_+$ and out of $M$ along $R_-$. Given a properly embedded surface
$F\subseteq M$ with a normal orientation and $\partial F$ transverse 
to $\gamma$, we can {\it decompose}
$(M,\gamma)$ along $F$ to obtain a new sutured manifold $(M|F,\gamma^\prime)$,
where the new sutures are constructed by an oriented sum of $\gamma$ and
$\partial F$, as in Figure 1. A degenerate (but important) example
is the decomposition $(E(K),\partial E(K)) \Rightarrow (E(\Sigma),\gamma)$
along $\Sigma$. A sequence of decompositions is a {\it sutured manifold
decomposition}. A sutured manifold is {\it taut} if 
it admits a sutured manifold decomposition ending with a disjoint union
of sutured manifolds of the form ($B^3,e$), where $e$ is the equatorial 
circle in the boundary of the 3-ball $B^3$. The  sequence of decomposing 
surfaces used in the decomposition is then called a
{\it taut sutured manifold heirarchy}. 
A fundamental theorem of Gabai [Ga1] states that
$(E(\Sigma),\gamma)$ admits a taut sutured manifold
heirarchy if and only if $\Sigma$ has minimal genus among all Seifert
surfaces for $K$. This gives, in principle, an effective way to compute the
genus of a knot or link. We will use this technology in Section 2 to 
compute the genera of our family of knots $K_n$.

\leavevmode

\epsfxsize=2.5in
\centerline{{\epsfbox{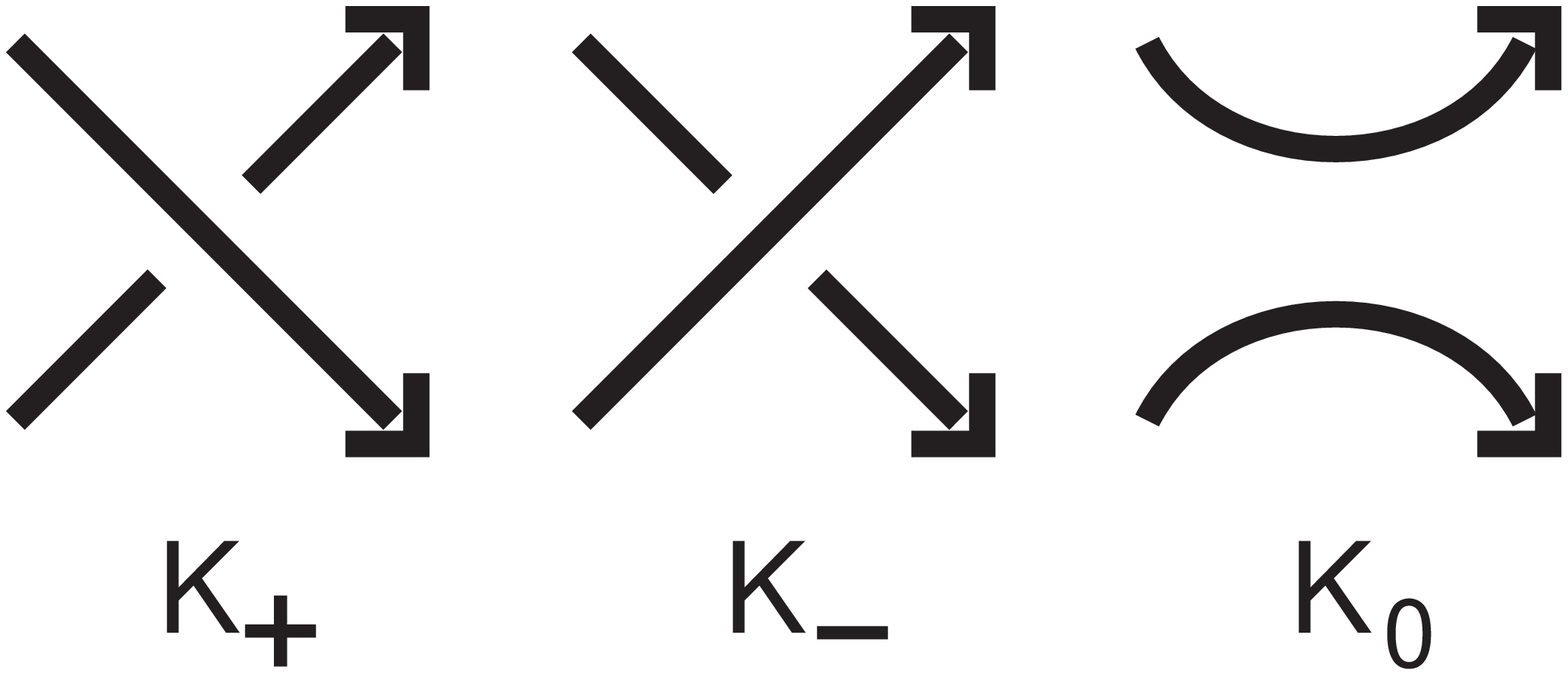}}}

\centerline{Figure 2: HOMFLY polynomial}

\medskip

The HOMFLY polynomial [FHLMOY] is a 2-variable Laurent polynomial defined 
for any oriented link, and may be thought of as the unique polynomial
$P_K(v,z)$,
defined on oriented link diagrams and invariant under the Reidemeister moves,
satisfying $P_{\text unknot}(v,z)=1$, and $v^{-1}P_{D_+} - vP_{D_-} = zP_{D_0}$,
where $D_+,D_-,D_0$ are diagrams for oriented links which all agree except 
at one crossing, where
they are given as in Figure 2. Here we are following the convention in
the naming of our variables found in Morton's paper. 

This skein relation gives an inductive method for computing the polynomial,
since one of $D_+,D_-$ will be ``closer'' than the other to the unlink,
in terms of unknotting number, while $D_0$ has fewer crossings.
It also allows any one of the polynomials to be computed from the other two.
We will use this relation in Section 3 to give upper bounds on the 
$z$-degree $M(K_n)$ of our family of knots $K_n$.

\leavevmode

\epsfxsize=4in
\centerline{{\epsfbox{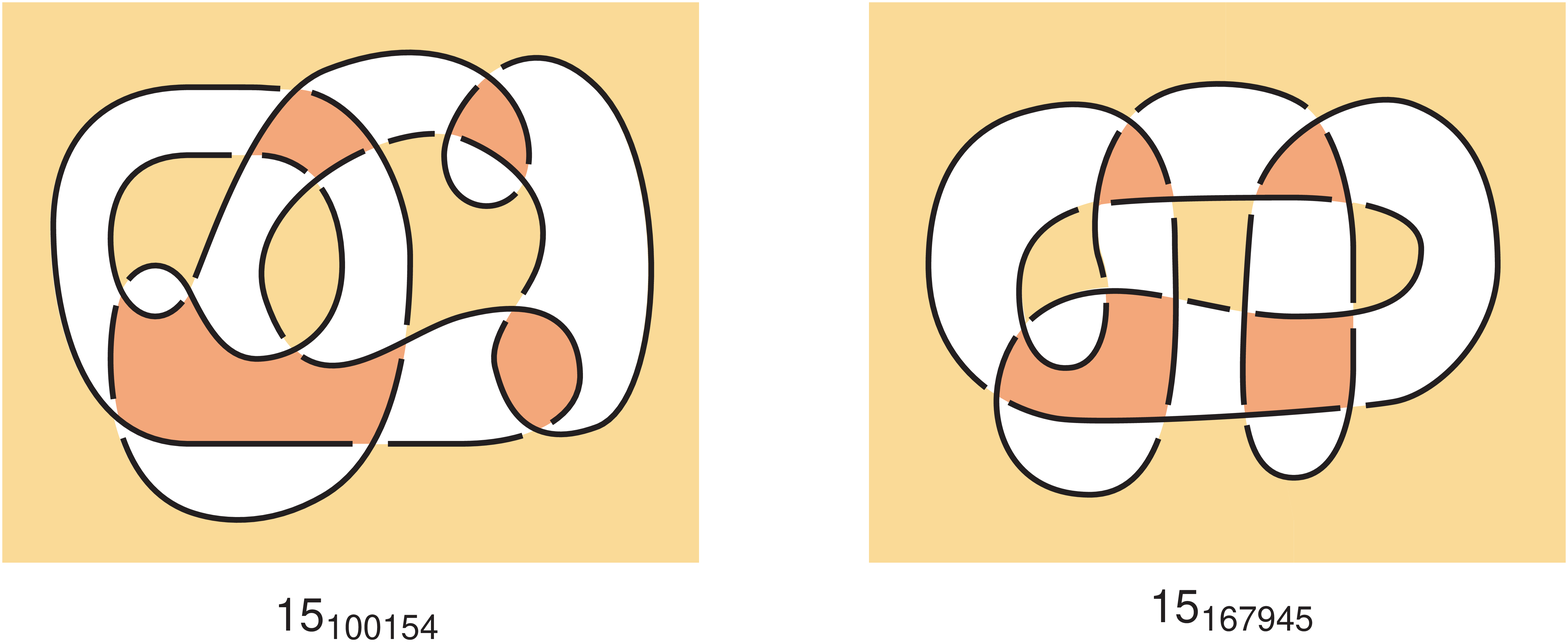}}}

\centerline{Figure 3: Stoimenow's examples}

\medskip

\heading{\S 2 \\ The canonical genus via sutured handlebodies}\endheading

The proof of our main theorem will be carried out in two steps.
In this section we will show that the knots $K_n$ described
in the theorem have the 
canonical genus claimed: $g_c(K_n) = g_c(K)+n$. In the next section
we will show that the HOMFLY polynomials of these knots have $z$-degree $M(K_n)$
strictly lower than this.

The basic idea behind Stoimenow's examples, and our argument, is that if
$g(K)$ is achieved by a canonical Seifert surface, then
this can be detected using sutured manifold theory. 
In particular, a canonical Seifert surface
$\Sigma$ has exterior $E(\Sigma)=E(K)|\Sigma$ a handlebody, which we
endow with a sutured manifold structure as above. If $(E(\Sigma),\gamma)$ 
admits a taut
sutured manifold heirarchy, then $\Sigma$ is a genus-minimizing, hence
$g_c$-minimizing, Seifert surface for $K$. For a sutured handlebody, the 
simplest heirarchy we can hope for is a {\it disk decomposition}, that
is, the intermediate decomposing surfaces are (compressing) disks for $E(\Sigma)$. 

\leavevmode

\epsfxsize=5in
\centerline{{\epsfbox{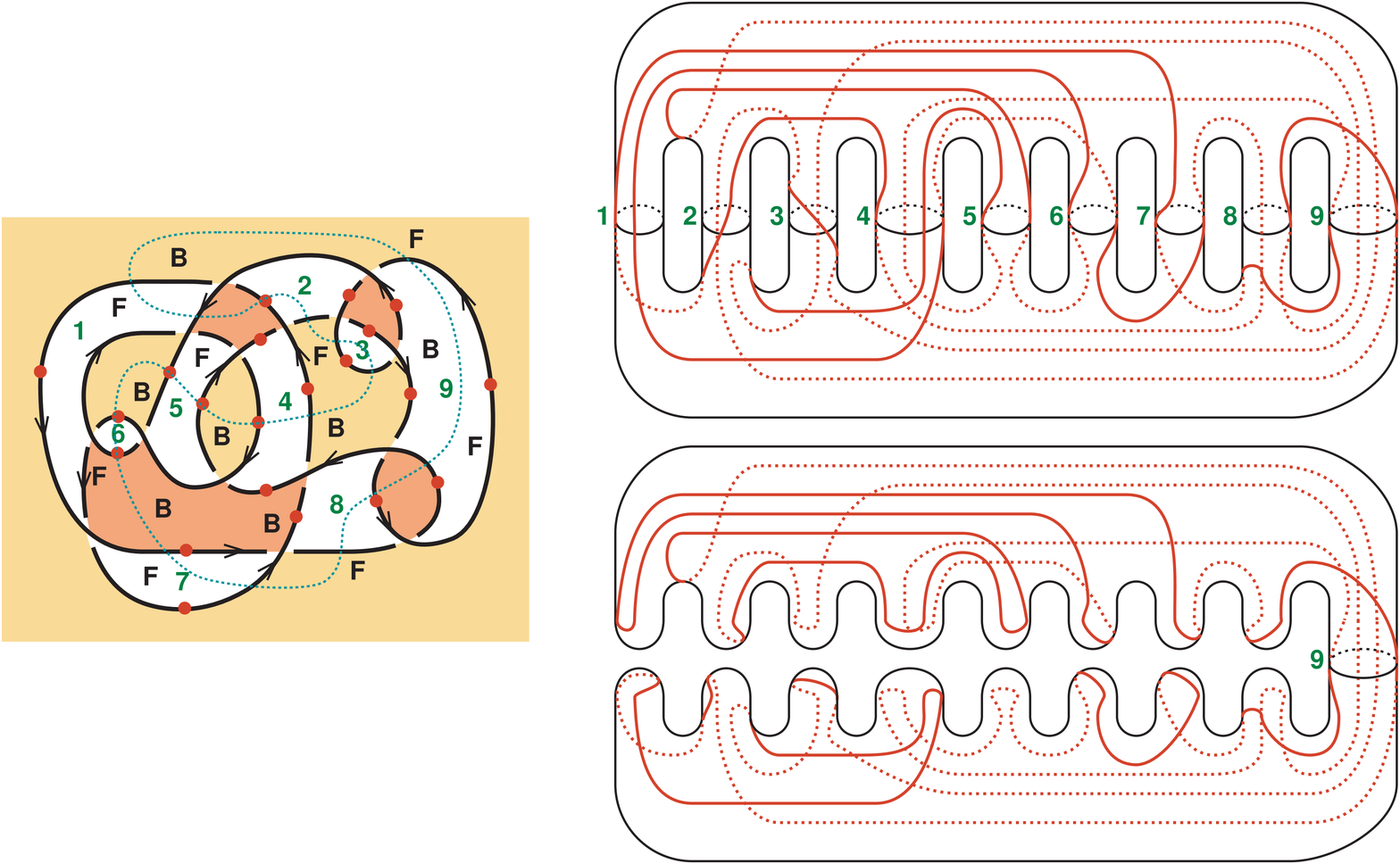}}}

\centerline{Figure 4: Sutured handlebodies}

\medskip

It is known that not all taut sutured handlebodies are disk decomposable [Go], and 
this phenomenon can even occur for Seifert surface exteriors [Br]. But in 
the case of Stoimenow's examples, a set of decomposing disks is readily
available. We show such a set of disks for the first example in 
Figure 4. The left-hand side of the figure describes a template 
for constructing
the sutured handlebody $(E(\Sigma),\gamma)$ as a standard handlebody, at the right; the
dotted line represents the intersection of a vertical plane of the 
paper running down the middle of the handlebody at the right.
``B'' marks the back of the figure at right, ``F'' the front. Looking 
into the paper at left is looking down from above at the right.
On the right-hand side, we illustrate the sutured manifold
resulting from decomposing along the horizontal compressing disks
labeled 1 through 8.
Since these disks cut $E(\Sigma)$ into a single 3-ball, and the suture under
decomposition becomes a single curve, the decomposition is taut,
so the canonical Seifert surface has minimal genus. A similar sequence of
images can be built for the other example given in Figure 3; 
we leave this for the interested reader.

The conditions described by the theorem, and the construction 
of $K_n$ from $K$, are as in Figure 6. In the discussion to follow, we assume 
for ease of exposition that the half-twisted band of the 
theorem joins a pair of Seifert disks that are not nested, i.e., the disks lie
in the same plane and are disjoint. This is not really an issue; every
canonical Seifert surface is isotopic to a checkerboard surface [Hi],
and the isotopy process may be assumed to leave fixed all crossings of our 
original diagram. The idea is that the isotopy to a checkerboard is
carried out inductively on the outermost disk of a nested collection,
as in Figure 5; this process adds, but does not delete, crossings to the 
underlying diagram. The crossing hypothesized by the theorem is therefore
still available to us in a checkerboard configuration.

\leavevmode

\epsfxsize=3in
\centerline{{\epsfbox{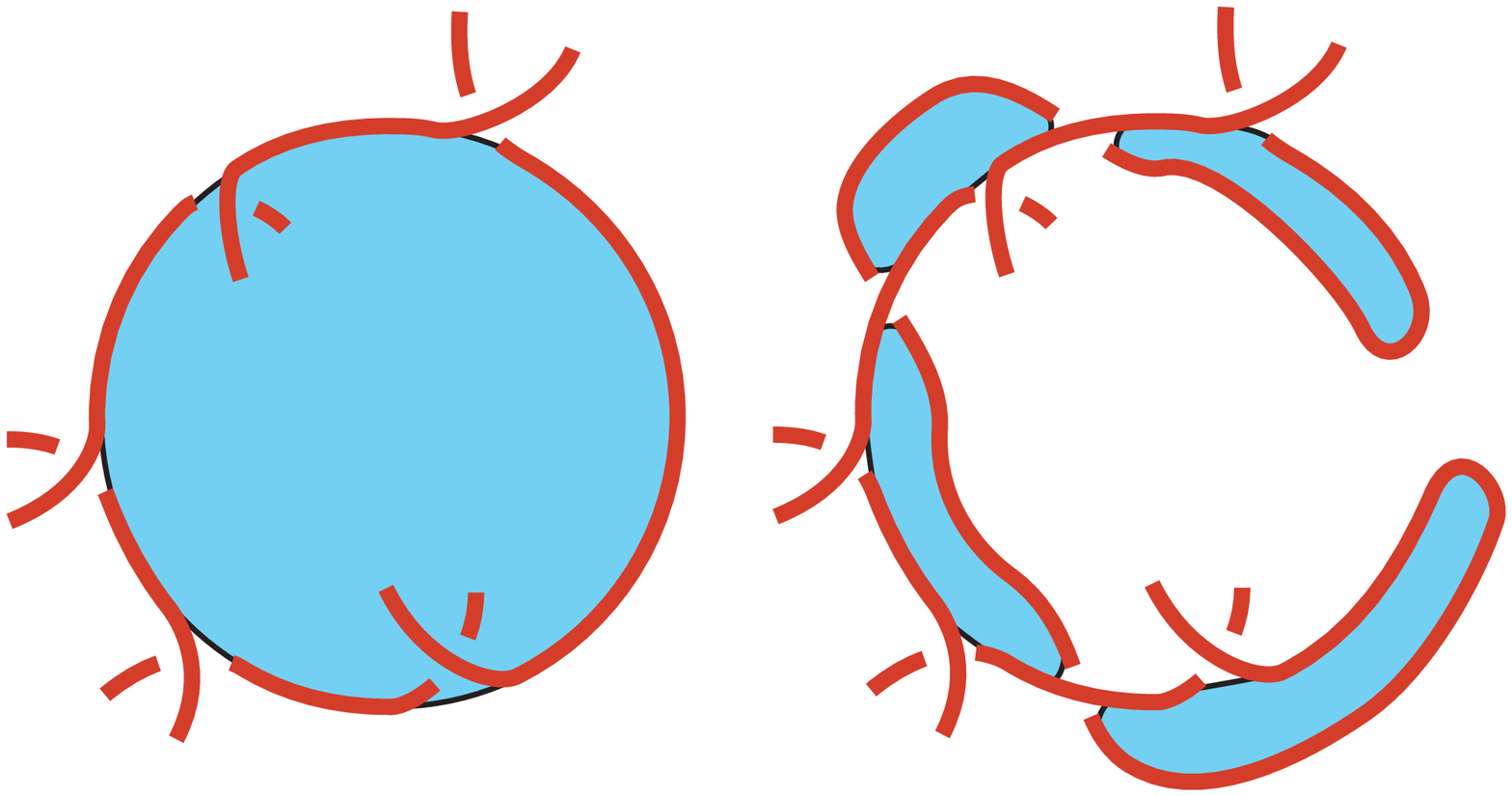}}}

\centerline{Figure 5: Making canonical surfaces checkerboard}

\medskip

To prove our claim that $g_c(K_n)=g_c(K)+n$, we show that the 
canonical Seifert surface $\Sigma_n$ for $K_n$ built from the 
diagram $D$ of the theorem admits a taut sutured 
manifold decomposition. Since each
additional pair of crossings introduced does not increase the number of Seifert
circles for $\Sigma_n$, the Euler characteristic decreases by two, so
the genus increases by one, each time.

\leavevmode

\epsfxsize=4.5in
\centerline{{\epsfbox{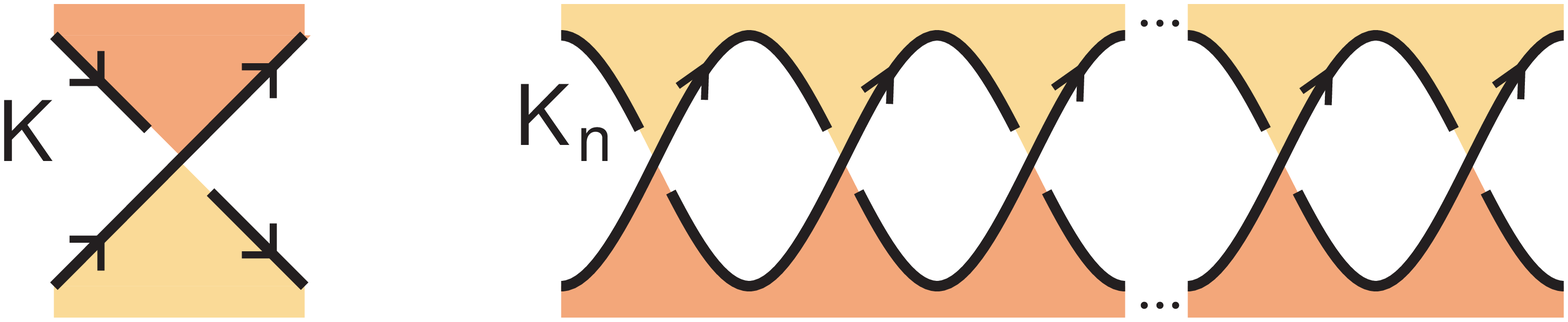}}}

\centerline{Figure 6: $K$ to $K_n$}

\medskip

Since $\Sigma$ is a genus-minimizing Seifert
surface for $K$, $(E(\Sigma),\gamma)$ admits a taut
sutured manifold decomposition. To show that $\Sigma_n$ is 
genus-minimizing for $K_n$, it suffices to show that there
is a sequence of decomposing surfaces taking $(E(\Sigma_n),\gamma_n)$
to $(E(\Sigma),\gamma)$; the taut sutured manifold heirarchy
for $(E(\Sigma),\gamma)$ can then be appended to this sequence
to give a taut heirarchy for $(E(\Sigma_n),\gamma_n)$. And to do
this, it suffices to find a sequence taking $(E(\Sigma_n),\gamma_n)$
to $(E(\Sigma_{n-1}),\gamma_{n-1})$. But this, in turn is
not difficult; the pair of compressing disks for the pair of 1-handles
dual to the pair of half-twisted bands added
to obtain $E(\Sigma_n)$ from $E(\Sigma_{n-1})$ provide the
necessary decomposing surface. On the level of sutured manifolds, this is
illustrated in Figure 7.  The basic idea is that since the disks ``look''
as if they belong in the exterior of the checkerboard surface of an 
alternating knot, the arguments
of [Ga2] ensure that the sutured manifold obtained by decomposing along
them is also the exterior of a Seifert surface; this surface is $\Sigma_{n-1}$.

By induction, therefore, we have that $E(\Sigma_n)$ admits a taut sutured manifold
heirarchy, and so, by induction, $g(K_n)=g(\Sigma_n)=g(\Sigma)+n$, as desired.

\leavevmode

\epsfxsize=4.5in
\centerline{{\epsfbox{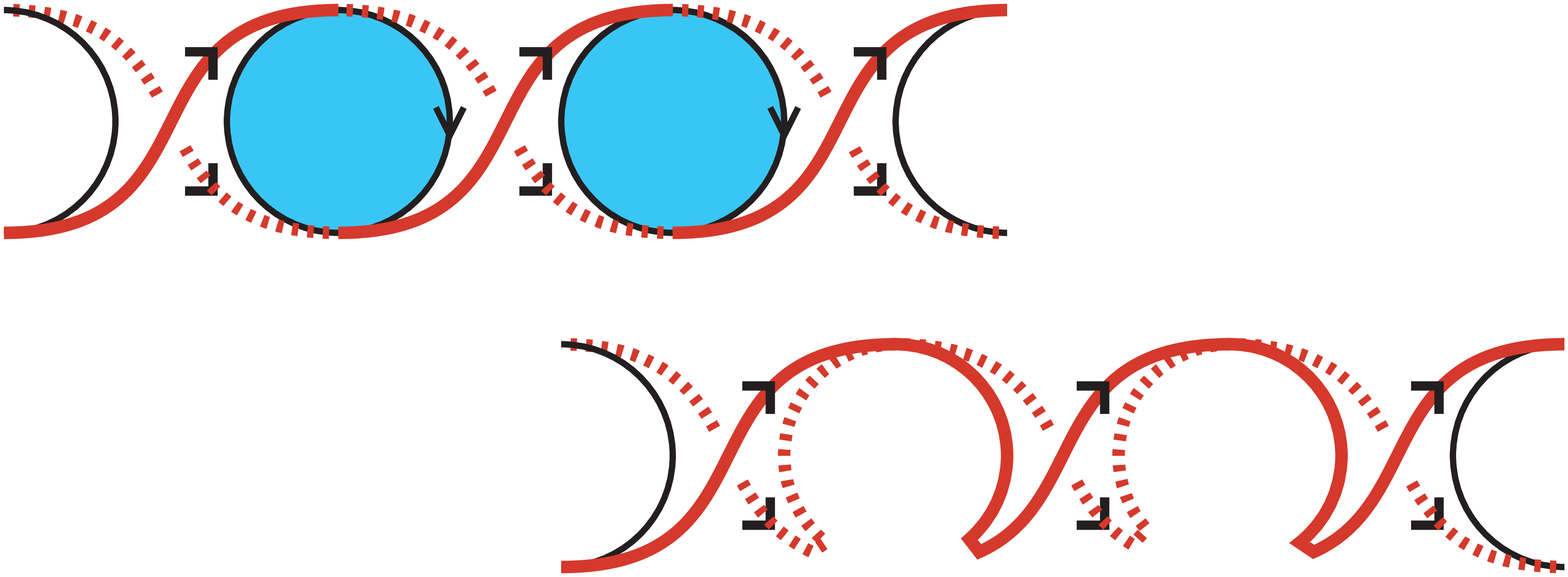}}}

\centerline{Figure 7: Heirarchies: the inductive step}

\medskip

\heading{\S 3 \\ The degree of the HOMFLY polynomial is too low}\endheading

In this section we finish the proof of our main theorem by showing that 
if $K$ is as in the theorem, then $M(K_n) < 2g_c(K) + 2n$ .
We then have $M(K_n) < 2g_c(K) + 2n = 2g_c(K_n)$,
as desired.
The needed inequality follows directly from the hypotheses of the theorem 
and an induction
argument using the skein relation satisfied by the HOMFLY polynomial.

\medskip

In what follows, we assume that the crossing given in the hypotheses of the theorem
is a positive crossing for $K$ and negative for $K^\prime$. A nearly identical
argument applies to the reverse situation. To avoid a possible confusion of notation, we will
denote by $L_n$ the link obtained by replacing the (positive) crossing $c$ with $n$
positive crossings joining the same two Seifert circles as $c$ does, in parallel. (So, in 
our theorem, $K=L_1$ and $K_n=L_{2n+1}$.) It then suffices to show that $M(L_n) < 2g_c(K)-1+n$,
since then $M(K_n) = M(L_{2n+1}) < 2g_c(K) -1+(2n+1) = 2g_c(K)+2n$.

Now by hypothesis, $g_c(K^\prime) < g_c(K)$. But by Morton's inequality, $M(K^\prime) 
\leq 2g_c(K^\prime)$, so $M(K^\prime) < 2g_c(K)$. 
From the skein equation for the HOMFLY polynomial, $v^{-1}P_{K_+} -vP_{K_-} = zP_{K_0}$, 
we can immediately conclude that, since the degree of the sum of two polynomials is no 
more than the maximum of their degrees,

\smallskip

\centerline{$M(K_0) \leq $max$\{M(K_+),M(K_-)\}-1$ 
\hskip10pt , \hskip10pt 
$M(K_+) \leq $max$\{M(K_-),M(K_0)+1\}$,}

\smallskip

\centerline{and $M(K_-) \leq $max$\{M(K_+),M(K_0)+1\}$.}

\smallskip

\noindent So using $K_+=K=L_1,K_-=K^\prime$, and $K_0=L_0$, we have

\smallskip

\centerline{$M(L_0)\leq$max$\{M(K),M(K^\prime)\}-1 < 2g_c(K)-1$ .}

\smallskip

\noindent This together with $M(L_1) = M(K) < 2g_c(K) =2g_c(K)-1+1$
give us the beginnings of our induction.

We proceed by complete induction. Assume that $n\geq 2$ and $M(L_j)<2g_c(K)-1+j$ for all $j<n$. Then
looking at any one of the $n$ parallel crossings $c$ of our diagram for $K_n$, we have 
$K_+=L_n$, $K_-=L_{n-2}$, and $K_0=L_{n-1}$. Then 

\smallskip

\noindent $M(K_+) = M(L_n) \leq $max$\{M(K_-),M(K_0)+1\} = 
$max$\{M(L_{n-2}),M(L_{n-1})+1\}$ 

\noindent $< $max$\{2g_c(K)-1+{n-2},(2g_c(K)-1+(n-1))+1\} = 2g_c(K)-1+n$, 

\smallskip

\noindent as desired.
Since the inequality is true for $n=0,1$, by complete induction, $M(L_n) < g_c(K) -1+n$ for all $n\geq 0$ .
This finishes the proof of Theorem 1. \blkbox

\leavevmode

\epsfxsize=4.5in
\centerline{{\epsfbox{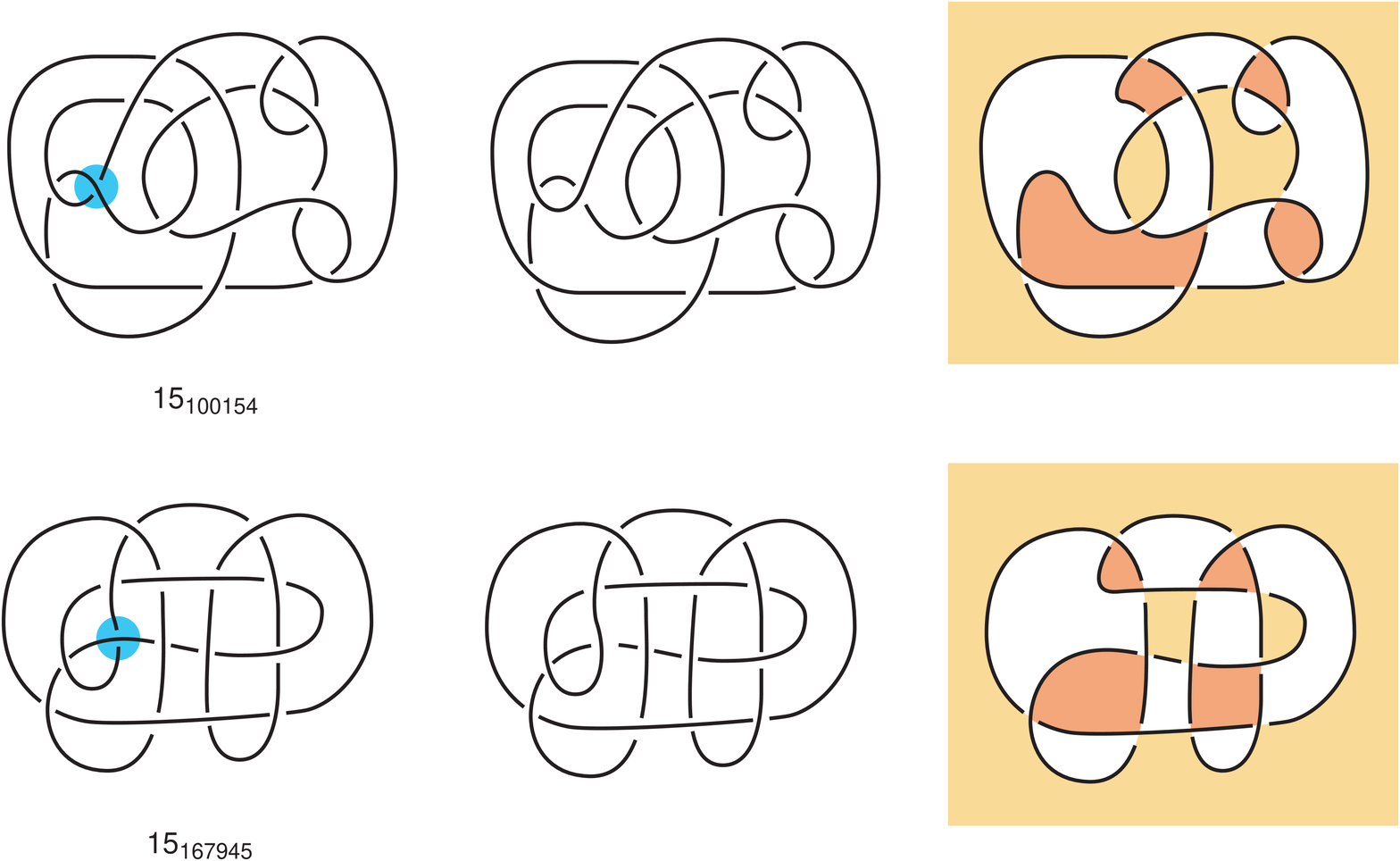}}}

\centerline{Figure 8: Explicit families of examples}

\medskip

We now prove Corollary 2 by demonstrating that Stoimenow's examples satisfy the hypotheses of the 
theorem. The sutured manifold calculations of section 2 verify that $g(K) = g_c(K) = 4$ for both knots.
A routine calculation shows that the HOMFLY polynomial of $K$ = $15_{100154}$ is
$(v^2+6v^{-2})z^6+(-v^4+4v^2+6-5v^{-2}+v^{-4})z^4+(-3v^4+4v^2+10-9v^{-2}+2v^{-4})z^2+(-2v^4+v^2+6-5v^{-2}+v^{-4})$ 
and for $K$ = $15_{167945}$ is 
$(v^2+v^{-2})z^6+(-v^4+2v^2+9-4v^{-2}-2v^{-4})z^4+(-2v^4+12-6v^{-2}-v^{-4}+v^{-6})z^2+(-v^4-v^2+6-3v^{-2})$ . 
In addition, each of the $g_c$-minimizing 
knot diagrams in Figure 3 has at least one crossing change which lowers the canonical genus;
see Figure 8 for canonical surfaces of genus 3. The relevant crossings join coplanar disks; 
both Seifert surfaces are 
checkerboards. Therefore all of the hypotheses of Theorem 1 are satisfied, and so each knot
gives rise to a family of knots $K_n$, for $n\geq 1$, with $g(K_n)=g_c(K_n)=n+3$ and $M(K_n)<2g_c(K_n)$.
This proves the corollary. \blkbox


\heading{\S 4 \\ Further considerations}\endheading

The main result of this paper was, in some sense, a by-product of our investigations
[BJ] into the canonical genus of Whitehead doubles of knots, motivated by the work of
Tripp [Tr] and Nakamura [Na]. (In what follows, we use the term ``the Whitehead double'' 
loosely here; in these investigations the resulting calculations are the same no matter 
how many twists the double has, unless $K$ is the unknot.) 
If $W(K)$ is the Whitehead double of the knot $K$, and $c(K)$ is the minimal crossing
number for $K$, then a direct construction of a canonical Seifert surface for $W(K)$
($K$ $\neq$ unknot) yields the inequality $g_c(W(K))\leq c(K)$ . Together with Morton's
inequality we then have 

\ctln{$M(W(K)) \leq 2g_c(W(K)) \leq 2c(K)$}

\noindent In the cases covered by the papers [Tr],[Na],[BJ], an induction argument
establishes that $M(W(K)) = 2c(K)$, proving that $g_c(W(K)) = c(K)$ and (as a by-product)
$M(W(K)) = 2g_c(W(K))$.

A posteriori, {\it this method of proof} requires that Morton's inequality be an equality.
Since we were aware of this prior to beginning our investigations in [BJ], we were led first 
to investigate the possible failure of Morton's ``equality'', which led us to the examples 
of Stoimenow and so to our own.
It is not too difficult to generate examples of knots $K$ for which
$M(W(K)) < 2c(K)$; the knots $K = 8_{19},8_{20}$ were the first that we found, with the
aid of Mathematica and the software package KnotTheory [KT]. (In fact, this inequality is 
strict for every non-alternating pretzel knot $K$ [BJ].) It follows that for these examples
{\it either} $M(W(K)) < 2g_c(W(K))$ {\it or} $g_c(W(K)) < c(K)$ . (Exactly one of these
is true, for each of the above two knots.) The question is: \underbar{which} of the two 
inequalities is strict? From the authors' point of view, strictness of the second one would
probably be the more interesting.

We are led to believe that Morton's inequality is the one which we would expect to be
an equality. This is supported by the known classes of knots for which it is true, mentioned
in the introduction. It is further supported by a result of Stoimenow [St3, Theorem 11.1];
he shows that, asymptotically, as the genus is held fixed but the number of crossings
goes to infinity, almost all canonical Seifert surfaces for a link diagram achieve 
the underlying link's genus, which in turn is equal to the degree of the Alexander polynomial.
Since, in terms of the HOMFLY polynomial, the Alexander polynomial 
$\Delta_K(t) = P_K(1,t^{1/2}-t^{-1/2})$,
the $z$-degree of the HOMFLY polynomial is at least 2$\cdot$deg$\Delta_K(t)$ = $2g_c(K)$,
in this case. Morton's inequality is therefore an equality.

It seems remarkable that the HOMFLY polynomial, which can be
computed from any projection of a knot $K$, can be so good at predicting the canonical
genus of $K$, and therefore, in some sense, good at predicting what the 
``simplest'' projection of
$K$ looks like (from the point of view of Seifert's algorithm). A better understanding of 
when the HOMFLY polynomial {\it fails} to do so, that is, when Morton's inequality is
strict, could help us better understand why it is usually so good at it. For example, can
one always find a skein tree diagram
in which every branch of the
tree exhibits an unexpected drop in degree as we progress to a collection of unlinks? 
Or do some knots, in the course of the calculation from a diagram, always have an 
unexpected collision 
of high-degree terms which fortuitously cancel one another out?

We should point out that we can, in principle, determine if Morton's inequality is an 
equality for any given knot $K$.
This is because we can recursively construct all knots having a canonical Seifert surface
up to a given genus $g$; all such surfaces are isotopic to checkerboard surfaces, and fall into
finitely-many twist-equivalent classes [St1]. HOMFLY polynomials within each class are related to
one another via the skein relation, and so we can quickly narrow our search down to
finitely-many knots having at most a given canonical genus and with the same HOMFLY 
polynomial as $K$. We can then use the solution 
to the homeomorphism problem for knot complements [He],[Wa] to check each candidate against $K$.
This is, of course, an extremely laborious process. But this might be worthwhile to carry out
for a single knot like $W(8_{19})$ or $W(8_{20})$, where we know that the canonical genus 
is at most 8 (by construction) and at least 7 (by Morton's inequality). Perhaps the best that
we could hope for is that no knot with canonical genus 7 has the same HOMFLY polynomial
as $W(8_{19})$ and/or $W(8_{20})$. Then we would not need to test for equivalence of knots;
Morton's inequality would necessarily be the inequality which is strict.

\medskip

{\bf Acknowledgements:} The first author wishes to thank the Department of Mathematics
of the City College of New York for their hospitality while a part of this work was 
carried out. The second author wishes to thank the Department of Mathematics of the
University of Nebraska - Lincoln and the Nebraska IMMERSE program for their 
hospitality and support.

\enddocument